\documentclass[11pt,fleqn]{article}
\usepackage{amsmath,amssymb}
\usepackage{wasysym}
\usepackage{graphicx}

\begin{document}
\sloppy

\newtheorem{axiom}{Axiom}[section]
\newtheorem{conjecture}[axiom]{Conjecture}
\newtheorem{corollary}[axiom]{Corollary}
\newtheorem{definition}[axiom]{Definition}
\newtheorem{example}[axiom]{Example}
\newtheorem{lemma}[axiom]{Lemma}
\newtheorem{observation}[axiom]{Observation}
\newtheorem{open}[axiom]{Problem}
\newtheorem{proposition}[axiom]{Proposition}
\newtheorem{theorem}[axiom]{Theorem}
\newcommand{\qed}{~~$\Box$}
\newcommand{\proof}{\emph{Proof.}\ \ }

\newcommand{\rz}{{\mathbb{R}}}
\newcommand{\zz}{{\mathbb{Z}}}
\newcommand{\cei}[1]{\lceil #1\rceil}
\newcommand{\flo}[1]{\left\lfloor #1\right\rfloor}

\newcommand{\tap}{\mbox{3AP}}
\newcommand{\mintap}{\mbox{min-3AP}}
\newcommand{\maxtap}{\mbox{max-3AP}} 
\newcommand{\ptap}{\mbox{planar-3AP}}
\newcommand{\tdm}{\mbox{3DM}}
\newcommand{\maxtdm}{\mbox{max-3DM}}
\newcommand{\pit}{\mbox{PIT}}

\newcommand{\aaa}{\alpha}
\newcommand{\bbb}{\beta}
\newcommand{\ccc}{\gamma}

\title{{\bf Geometric versions of the 3-dimensional assignment problem under general norms}}
\author{\sc Ante {\'C}usti{\'c}
\thanks{{\tt custic@opt.math.tugraz.at}.
Institut f\"ur Optimierung und Diskrete Mathematik, TU Graz, Steyrergasse 30, A-8010 Graz, Austria}
\and
\sc Bettina Klinz
\thanks{{\tt klinz@opt.math.tugraz.at}.
Institut f\"ur Optimierung und Diskrete Mathematik, TU Graz, Steyrergasse 30, A-8010 Graz, Austria}
\and
\sc Gerhard J.\ Woeginger
\thanks{{\tt gwoegi@win.tue.nl}.  Department of Mathematics and Computer Science,
TU Eindhoven, P.O.\ Box 513, 5600 MB Eindhoven, Netherlands}}
\date{}
\maketitle

\begin{abstract}
We discuss the computational complexity of special cases of the 
3-dimensional (axial) assignment problem 
where the elements are points in a Cartesian space and where the cost coefficients are the perimeters 
of the corresponding triangles measured according to a certain norm.
(All our results also carry over to the corresponding special cases of the 3-dimensional 
matching problem.)

The minimization version is NP-hard for every norm, even if the underlying Cartesian space is 2-dimensional.
The maximization version is polynomially solvable, if the dimension of the Cartesian space is fixed 
and if the considered norm has a polyhedral unit ball. 
If the dimension of the Cartesian space is part of the input, the maximization version is NP-hard for
every $L_p$ norm; in particular the problem is NP-hard for the Manhattan norm $L_1$ and the Maximum 
norm $L_{\infty}$ which both have polyhedral unit balls.

\end{abstract}

\section{Introduction}
\nopagebreak
The 3-dimensional (axial) assignment problem (\tap) is an important 
and well-studied problem in combinatorial optimization.
An instance of the {\tap} consists of three sets $X$, $Y$, $Z$ with $|X|=|Y|=|Z|=n$, and a cost function 
$c:X\times Y\times Z\to\rz$.
The goal is to find a set of $n$ triples in $X\times Y\times Z$ that cover every element in $X\cup Y\cup Z$
exactly once, such that the sum of the costs of these triples is minimized.
In the closely related maximization version {\maxtap} of the {\tap}, this sum is to be maximized.
The book \cite{Burkard-book} by Burkard, Dell'Amico \& Martello contains a wealth of information on the
{\tap} and other assignment problems.

A prominent special case of the {\tap} is centered around some metric space $(S,d)$ where $S$ is a set and 
where $d$ is a distance function on $S\times S$ (that hence is symmetric, non-negative, and satisfies 
the triangle inequality).
The elements in $X\cup Y\cup Z$ are points in $S$, and the cost $c(x,y,z)$ of a triple 
$(x,y,z)\in X\times Y\times Z$ is given by
\begin{equation}
\label{eq:perimeter}
c(x,y,z) ~=~ d(x,y)+d(y,z)+d(z,x).
\end{equation}
Costs of this type are called \emph{perimeter costs}; intuitively speaking, they measure the perimeter
of the triangle determined by three points $x,y,z$ in the metric space.

The {\tap} is well-known to be NP-hard; see for instance Karp~\cite{Karp1972} or Garey \& Johnson \cite{GaJo1979}.
Spieksma \& Woeginger \cite{SpWo1996} establish NP-hardness of the special case of perimeter costs (\ref{eq:perimeter})
where the underlying metric space is the two-dimensional Euclidean plane with standard Euclidean distances.
Polyakovskiy, Spieksma \& Woeginger \cite{PoSpWo2013} show that {\tap} and {\maxtap} with perimeter costs are 
polynomially solvable, if the underlying metric space satisfies the so-called Kalmanson conditions; their results
cover convex Euclidean point sets and tree metric spaces as special cases.
Crama \& Spieksma \cite{CrSp1992} design a polynomial time approximation algorithm with worst case guarantee $4/3$
for the {\tap} with perimeter costs; their approach works for arbitrary metric spaces without imposing any
additional structural constraints.
Burkard, Rudolf \& Woeginger \cite{BuRuWo1996} exhibit a polynomially solvable special case of the {\maxtap} 
where the costs are decomposable and products of certain parameters.

\paragraph{Results of this paper.}
We study {\tap} and {\maxtap} with perimeter costs in Cartesian spaces under arbitrary distance functions.
On the negative side, we derive NP-hardness results that contain and generalize the known results from the 
literature for the standard Euclidean distances.
On the positive side, we derive polynomial time algorithms for certain special cases of {\maxtap} where 
the distances are defined via norms with polyhedral unit balls.
Our main results are the following:
\begin{itemize}
\itemsep=0.2ex
\item[(A)]
Problem {\maxtap} is polynomially solvable, if the dimension of the underlying Cartesian space is a fixed
constant and if the underlying norm has a polyhedral unit ball.
\item[(B)]
Problem {\maxtap} is NP-hard, if the dimension of the underlying Cartesian space is part of the input and
if the underlying norm is any fixed $L_p$ norm. 
This hardness result in particular holds for the Manhattan norm $L_1$ and the Maximum norm $L_{\infty}$ 
which both have polyhedral unit balls.
\item[(C)]
Finally, the minimization problem {\tap} is NP-hard for any fixed norm, even if the underlying Cartesian 
space is $2$-dimensional.
\end{itemize}
Result (A) heavily builds on the machinery developed by Barvinok, Fekete, Johnson, Tamir, 
Woeginger \& Woodroofe \cite{Barvinok2003} for the Travelling Salesman 
Problem (TSP).
Also the TSP is polynomially solvable, if the cities are points in some Cartesian space of fixed
dimension and if the distances are defined via norms with polyhedral unit balls.
While the framework for our result (A) is taken from \cite{Barvinok2003}, the technical details and
the combinatorial features are very different and require a number of new ideas.
Result (B) is done by a routine NP-hardness reduction from a closely related NP-hard graph problem.
Result (C) builds on the NP-hardness reductions of Spieksma \& Woeginger \cite{SpWo1996} and
Pferschy, Rudolf \& Woeginger \cite{PfRuWo1994} for Euclidean distances.
In the Euclidean case, one may use Pythagorean triangles as simple building blocks to control the 
distances between points and to ensure rational coordinates that can be processed by a Turing machine.
In the general case (C), it is much more tedious to prove the existence of the corresponding building blocks.

\paragraph{Organization of this paper.}
Section~\ref{sec:preliminaries} summarizes some standard geometric definitions around distances, norms
and unit balls.
Result (A) for the {\maxtap} is derived in two steps.
First Section~\ref{sec:tunnel} derives an auxiliary result on the {\maxtap} under so-called tunneling 
distances, and then Section~\ref{sec:polynorm} establishes that {\maxtap} under polyhedral norms is
a special case of the tunneling case.
Section~\ref{sec:lp-norm} contains the proof of result (B).
Section~\ref{sec:lattice} constructs certain lattices with certain useful properties; these lattices
are then used in Section~\ref{sec:min} to prove the NP-hardness result (C).
Section~\ref{sec:3dm} translates our results (A), (B) and (C) into corresponding results for the 
maximization version and the minimization version of the 3-dimensional matching problem.
Finally, Section~\ref{sec:conclusions} concludes the paper with a short discussion and some open questions.

\bigskip
\section{Technical preliminaries}
\label{sec:preliminaries}
\begin{figure}[bht]
\bigskip
\centerline{\includegraphics[height=3.3cm]{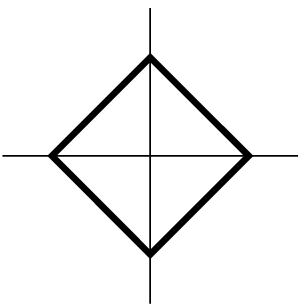} \hspace{0.8cm}
            \includegraphics[height=3.3cm]{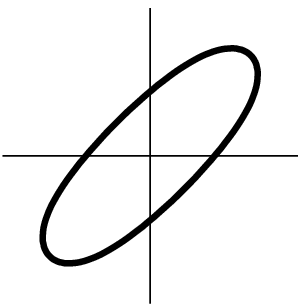} \hspace{0.8cm}
            \includegraphics[height=3.3cm]{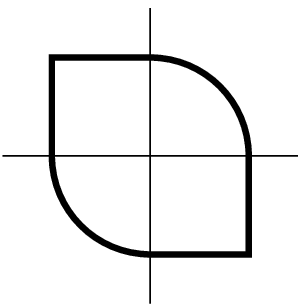}}
\caption{Three examples of unit balls for an $L_R$ norm in $\rz^2$}
\label{fig:ball}
\end{figure}

\nopagebreak
Let $R$ denote a compact and convex subset of the $s$-dimensional Cartesian space $\rz^s$ that has non-empty
interior and that is centrally symmetric with respect to the origin.
The corresponding norm $L_R$ with unit ball $R$ determines for any two points $x,y\in\rz^s$ a distance 
$d_R(x,y)$ in the following way.
First translate the space so that one of the two points (say point $x$) lies in the origin. 
Then determine the unique scaling factor $\lambda$ by which one must rescale the unit ball $R$ (shrinking 
for $\lambda<1$, expanding for $\lambda>1$), such that the other point (point $y$, in our case) lies on 
its boundary.
The distance is then given by $d_R(x,y)=\lambda$.
Note that since $R$ is centrally symmetric, it does not matter whether we choose point $x$ or point $y$ 
for the origin.
See Figure~\ref{fig:ball} for an illustration.

The most popular norms for $\rz^s$ are the Manhattan norm, the Euclidean norm, and the Maximum norm.
These three norms are special cases of the well-known $L_p$ norm, respectively for $p=1$, for $p=2$, 
and for $p=\infty$.
We recall that for $1\le p<\infty$, the $L_p$ distance between two points $x=(x_1,\ldots,x_s)$ and 
$y=(y_1,\ldots,y_s)$ in $s$-dimensional space is given by
\begin{equation}
\label{eq:Lp-norm}
d(x,y) ~=~ \left( \sum_{i=1}^s |x_i-y_i|^p \right)^{1/p}.
\end{equation}
For $p=\infty$, the corresponding distance under the Maximum norm $L_{\infty}$ is given by
\begin{equation}
\label{eq:L-infty-norm}
d(x,y) ~=~ \max{}_{i=1}^s |x_i-y_i|.
\end{equation}

\bigskip
\section{The maximization problem under tunneling distances}
\label{sec:tunnel}
\nopagebreak
In this section we consider a variant of the {\maxtap} with perimeter costs that will be useful in 
Section~\ref{sec:polynorm} of the paper.
The distances between the elements of $X\cup Y\cup Z$ are specified with the help of a system of $k\ge2$ 
so-called \emph{tunnels} $t_1,\ldots,t_k$; we stress that throughout this section the number $k$ of
tunnels is a constant that does not depend on the input.
Each tunnel acts as a bidirectional passage with a front entry and a back entry.
For every element $x\in X\cup Y\cup Z$ and every tunnel $t$, we denote 
by $F(x,t)$ the distance between $x$ and the front entry of $t$ and
by $B(x,t)$ the distance between $x$ and the back  entry of $t$.
Intuitively speaking, the only way of moving from $x$ to $y$ is to first move from $x$ to some tunnel, 
then to traverse the tunnel in either direction (either from front entry to back entry, or from back 
entry to front entry), and finally to move from the other end of the tunnel to $y$.
The \emph{tunneling distance} between two elements $x$ and $y$ in $X\cup Y\cup Z$ is then given by
\begin{equation}
\label{eq:tunnel.1}
d(x,y) ~=~ \max\left\{ F(x,t_i)+B(y,t_i),~ B(x,t_i)+F(y,t_i):~ 1\le i\le k\right\}.
\end{equation}
(We note in passing that the lengths of the tunnels do not play any role in this formula, 
as these lengths can easily be encoded in the values $F(x,t)$ and $B(x,t)$.)
 
We construct an undirected, edge-labeled, bipartite multigraph $G$ whose vertex set are the 
elements of $X\cup Y\cup Z$ together with the tunnels $t_1,\ldots,t_k$.
Between any element $x$ of $X\cup Y\cup Z$ and any tunnel $t$ there are four 
edges, two of which are
labeled $B$ and have cost $B(x,t)$, whereas the other two are labeled $F$ and 
have cost $F(x,t)$.

A \emph{six-cycle} in $G$ is a closed walk $x-t_i-y-t_j-z-t_{\ell}-x$ 
with $(x,y,z)\in X\times Y\times Z$ 
and three (not necessarily distinct) tunnels $t_i,t_j,t_{\ell}$. 
The six-cycle is \emph{legal}, if the labels of the two edges incident to $t_i$ are distinct, if the 
labels of the two edges incident to $t_j$ are distinct, and if the labels of the two edges incident 
to $t_{\ell}$ are distinct.
We stress that we do not require tunnel $t_i$ to be the maximizer of the expression for $x$ and $y$ 
in the right hand side of (\ref{eq:tunnel.1}), nor that $t_j$ and $t_{\ell}$ 
are the maximizers for the
corresponding expressions for $y$ and $z$, respectively for $z$ and $x$.

A \emph{legal} set $C$ of six-cycles consists of $n$ legal six-cycles in $G$ that cover every vertex
of $X\cup Y\cup Z$ exactly once.
We define $G[C]$ as the subgraph of $G$ that is induced by the $6n$ edges in $C$.
Then we coarsen the subgraph $G[C]$ by anonymizing the identities of the vertices in $X\cup Y\cup Z$:
every vertex in $X$ is simply labeled $X$,
every vertex in $Y$ is labeled $Y$, and
every vertex in $Z$ is labeled $Z$.
The resulting anonymized graph $G^*[C]$ is called an \emph{outline} of $C$ and $G[C]$.

\begin{lemma}
\label{le:tunnel.1}
The optimal objective value of the considered {\maxtap} instance coincides with the largest cost 
taken over all subgraphs $G[C]$ of $G$ with a legal set $C$ of six-cycles.
\end{lemma}
\proof
Let $C$ be an arbitrary legal set of six-cycles.
Every six-cycle $x-t_i-y-t_j-z-t_{\ell}-x$ in $C$ yields a corresponding triple $(x,y,z)$ in $X\times Y\times Z$.
The cost $c(x,y,z)$ of triple $(x,y,z)$ may be computed according to (\ref{eq:tunnel.1}), by replacing the
three tunnels $t_i,t_j,t_{\ell}$ by three other tunnels that maximize the value.
Hence, the cost of the triple is an upper bound on the cost of the six-cycle, and the cost of all $n$ 
corresponding triples is an upper bound on the cost of $G[C]$.
This shows that the optimal objective {\maxtap} value is an upper bound on the cost of every subgraph $G[C]$.

Next, consider a set $T$ of $n$ triples in $X\times Y\times Z$ that constitutes an optimal solution for the 
{\maxtap} instance.
We translate every triple $(x,y,z)\in T$ into a legal six-cycle: we let $t_i$ (respectively, $t_j$ and $t_{\ell}$)
denote the tunnel that maximizes the expression (\ref{eq:tunnel.1}) for $x$ and $y$ (respectively, for $y$
and $z$ and for $z$ and $x$), and we choose the labels $B$ and $F$ appropriately in the obvious way.
For the resulting legal set $C_T$, the cost of $G[C_T]$ coincides with the optimal objective {\maxtap} value.
\qed 

\begin{lemma}
\label{le:tunnel.2}
Let $G^*$ be a given outline.
Then one can compute in polynomial time $O(n^3)$ the largest cost of all the induced subgraphs $G[C]$ 
(with a legal set $C$ of six-cycles), whose outline $G^*[C]$ coincides with $G^*$.
\end{lemma}
\proof
The problem boils down to assigning the elements of $X$ (respectively, of $Y$ and $Z$) to the $n$ 
vertices in $G^*$ that are labeled $X$ (respectively, labeled $Y$ and $Z$).
The cost of assigning an element $x\in X$ to some vertex $v$ only depends on $x$ and on the 
two edges incident to $v$ in $G^*$.
Hence, we are dealing with a classical two-dimensional assignment problem which can be solved in 
polynomial time $O(n^3)$; see for instance Burkard, Dell'Amico \& Martello \cite{Burkard-book}.
\qed

\begin{lemma}
\label{le:tunnel.3}
There exist only $O(n^{8k^3})$ distinct outlines $G^*$ for graph $G$, and they can all be enumerated 
in polynomial time.
\end{lemma}
\proof
After anonymizing the identities of the vertices in $X\cup Y\cup Z$, a legal six-cycle is determined by
the three tunnels $t_i,t_j,t_{\ell}$ and the labels of its first, third, and fifth edge.
Hence there remain only $8k^3$ combinatorially different legal six-cycles, and each of them may be
used at most $n$ times in any outline.
\qed

\bigskip
The three lemmas suggest the following approach to {\maxtap} under tunneling distances:
enumerate all possible outlines in polynomial time according to Lemma~\ref{le:tunnel.3}, and for each such 
outline compute the maximum possible cost of a corresponding induced subgraph according to Lemma~\ref{le:tunnel.2}.
Return the largest cost over all outlines, which by Lemma~\ref{le:tunnel.1} coincides with the optimal
objective value of the {\maxtap} instance.

\begin{theorem}
\label{th:tunnel.main}
Problem {\maxtap} with perimeter costs under tunneling distances can be solved within a time complexity that 
depends polynomially on the instance size $n$ (and exponentially on the number
$k$ of tunnels).
\end{theorem}
\proof
The above approach computes the optimal objective value and the corresponding graphs $G^*[C]$ and $G[C]$,
but does neither yield the corresponding optimal solution $T\subset X\times Y\times Z$ for the {\maxtap} 
instance nor the underlying legal set $C_T$ of six-cycles.
We briefly sketch how these objects can also be determined in polynomial time.
The set $C_T$ can be determined in polynomial time by invoking Lenstra's algorithm \cite{Lenstra1983} for
integer programming in constant dimension.
For each of the $8k^3$ combinatorially different legal six-cycles, we introduce a corresponding integer
variable that counts the number of occurrences of this cycle in $C_T$.
The constraints in the integer program enforce that $G[C_T]$ coincides with $G[C]$.
And once we have found $C_T$ through the integer program, it is straightforward to identify the optimal
solution $T$ (as outlined in the proof of Lemma~\ref{le:tunnel.1}).
\qed

\bigskip
\section{The maximization problem under polyhedral norms}
\label{sec:polynorm}
\nopagebreak
Throughout this section, we consider the $s$-dimensional Cartesian space $\rz^s$ endowed with some 
fixed norm with polyhedral unit ball $R$.
We investigate the special case of {\maxtap} with perimeter costs where the elements in $X\cup Y\cup Z$ are 
points in $\rz^s$ and where the distances are measured according to $d_R$.
We stress that both the dimension $s$ of the underlying space and the number of faces of the unit ball $R$ 
are constants that do not depend on the input.

The unit ball $R$ is a polytope with $2k$ faces that is centrally symmetric with respect to the origin.
Then for certain vectors $h_1,\ldots,h_k\in\rz^s$, this polytope $R$ can be written as the intersection 
of a collection of half-spaces:
\begin{equation}
\label{eq:polynorm.1}
R ~=~ 
\left( \bigcap_{i=1}^k \{x:~ h_i\cdot x\le 1\} \right) \,\cap\, 
\left( \bigcap_{i=1}^k \{x:~ h_i\cdot x\ge-1\} \right)
\end{equation}
As an example, for the Manhattan norm in $\rz^2$ the corresponding vectors are
$h_1=(1,1)$ and $h_2=(-1,1)$, and for the Maximum norm in $\rz^2$ the corresponding vectors are
$h_1=(1,0)$ and $h_2=(0,1)$.
The distance $d_R(x,y)$ between two points $x,y\in\rz^s$ may then be written as
\begin{eqnarray}
d_R(x,y) 
&=&  \max \left\{ \left|h_i\cdot(x-y)\right|:~ 1\le i\le s \right\}     \nonumber \\[0.5ex]
&=&  \max \left\{ h_i\cdot(x-y),~ h_i\cdot(y-x) :~ 1\le i\le s \right\} \nonumber \\[0.5ex]
&=&  \max \left\{ h_i\cdot x-h_i\cdot y,~ -h_i\cdot x+h_i\cdot y :~ 1\le i\le s \right\} \label{eq:polynorm.2}
\end{eqnarray}
We model a {\maxtap} instance under a polyhedral norm as a special instance of {\maxtap} under tunneling 
distances as discussed in Section~\ref{sec:tunnel}.
The $k$ vectors $h_1,\ldots,h_k$ serve as tunnels, and we set $F(x,h_i)=x\cdot h_i$ and $B(x,h_i)=-x\cdot h_i$.
With this choice, the polyhedral distance $d_R(x,y)$ between two points $x$ and $y$ in $X\cup Y\cup Z$ 
in (\ref{eq:polynorm.2}) coincides with the tunneling distance given in (\ref{eq:tunnel.1}).
Hence Theorem~\ref{th:tunnel.main} yields the following.

\begin{theorem}
\label{th:polynorm.main}
For any polyhedral norm $L_R$ with unit ball $R$ in $s$-dimensional space $\rz^s$, problem {\maxtap} with 
perimeter costs measured according to $L_R$ can be solved within a time complexity that depends polynomially 
on the instance size $n$ (and exponentially on the number $k$ of facets of the polyhedral unit ball).
\qed
\end{theorem}

Theorem~\ref{th:polynorm.main} also implies the existence of a polynomial time approximation scheme (PTAS)
for {\maxtap} under any arbitrary norm with a not necessarily polyhedral unit ball $R$. 
One simply approximates the unit ball $R$ by a polyhedral unit ball. 
Since the dimension $s$ of the underlying space and the ball $R$ are fixed, one may choose a fixed
polyhedral approximation of the ball that approximates the distances between any two points within a
factor $1\pm\varepsilon$.
(This trick of approximating the unit ball by a polyhedral unit ball is essentially due to 
Barvinok \cite{Barvinok1996} who applied it to the maximum Travelling Salesman 
Problem.)

\begin{theorem}
\label{th:polynorm.ptas}
For any fixed (not necessarily polyhedral) norm $L_R$ with unit ball $R$ in $s$-dimensional space $\rz^s$, 
problem {\maxtap} with perimeter costs measured according to $L_R$ possesses a PTAS.
\qed
\end{theorem}

\bigskip
\section{The maximization problem in non-fixed dimension}
\label{sec:lp-norm}
\nopagebreak
The polynomial time results for {\maxtap} in the preceding section assumed that the dimension $s$ of 
the underlying Cartesian space $\rz^s$ as well as the number of faces of the underlying unit ball 
are constants that do not depend on the input.
In this section we discuss problem {\maxtap} with perimeter costs measured according to a standard
$L_p$ norm (with $1\le p\le\infty$) when the dimension $s$ is not fixed, but 
part of the input.
Our reductions are from the following variant of Partition into Triangles.

\begin{quote}
Problem: Partition into Triangles ({\pit})
\\[1.0ex]
Instance:
A $6$-regular, tripartite graph $G=(V,E)$ with tripartition $V=V_1\cup V_2\cup V_3$, where
$|V_1|=|V_2|=|V_3|=q$.
\\[1.0ex]
Question:
Does there exist a set $T$ of $q$ triples in $V_1\times V_2\times V_3$ such that every vertex 
in $V$ occurs in exactly one triple and such that every triple induces a triangle in $G$?
\end{quote}
We have not been able to locate an NP-hardness proof of {\pit} on $6$-regular tripartite graphs
in the literature (though we strongly expect that this result has been observed before).
For instance Van Rooij, Van Kooten Niekerk \& Bodlaender \cite{VanRooij2013} establish NP-hardness 
for $4$-regular graphs, but their graphs are not tripartite.

\begin{proposition}
\label{pr:pit}
Problem {\pit} on $6$-regular tripartite graphs is NP-complete.
\end{proposition}
\proof
The argument is routine, and we only sketch the main ideas.
The NP-hardness proof on pages 68 and 69 of Garey \& Johnson \cite{GaJo1979} 
for Partition into 
Triangles is a reduction from the Exact Cover By 3-Sets problem.
We perform essentially the same reduction, but start it from another NP-hard feasibility version 
of the 3-dimensional assignment problem with bounded occurrence of elements (Instance: three 
sets $X$, $Y$, $Z$ with $|X|=|Y|=|Z|=q$, and a set $T\subseteq X\times Y\times Z$ of triples such 
that every element of $X\cup Y\cup Z$ occurs in at most three triples of $T$. 
Question: Does there exist a subset $T^*$ of $q$ triples in $T$ such that each element of 
$X\cup Y\cup Z$ is contained in precisely one triple of $T^*$?).
Then the resulting graph $G$ is tripartite and all vertex degrees lie in $\{3,4,5,6\}$.

Hence it remains to make the graph $6$-regular.
This can be reached by various gadget constructions.
We sketch a particularly simple approach that increases the minimum degree of $G$ by~$1$, while
keeping the maximum degree unchanged.
Take the graph $G=(V,E)$, and construct a copy $G'=(V',E')$ of it (so that for every $v\in V$ there
is a corresponding copy $v'\in V'$, and there is an edge $[u,v]\in E$ if and only if there is an
edge $[u',v']\in E'$).
Define a new graph on the vertex set $V\cup V'$, and all edges in $E\cup E'$, and furthermore an
additional edge between $v$ and $v'$ whenever vertex $v$ has degree in $\{3,4,5\}$.
The new graph is still tripartite, and it has a partition into triangles if and only if the old
graph allows a partition into triangles (note that the additional edges $[v,v']$ do not occur in 
any triangle, and hence are irrelevant for partitions into triangles).
If we repeat this construction two more times, the resulting graph will be $6$-regular and tripartite.
\qed

\bigskip
The following two lemmas establish NP-hardness of {\maxtap} with perimeter costs for all
values $p$ with $1\le p\le\infty$.

\begin{lemma}
\label{lp-norm.1}
For any fixed $p$ with $1\le p<\infty$, problem {\maxtap} with perimeter costs measured according 
to the $L_p$ norm is NP-hard.
\end{lemma}
\proof
We consider an arbitrary instance $G=(V,E)$ of {\pit} with $|V|=3q$, and we construct the following 
instance of {\maxtap} with perimeter costs from it.
For every vertex $v$ in part $V_1$ (respectively, part $V_2$ and part $V_3$), we create a corresponding 
point $P(v)$ that belongs to the set $X$ (respectively, set $Y$ and set $Z$).
We choose the dimension $s=\binom{3q}{2}$, and we make every coordinate correspond to one $2$-element set 
of vertices in $V$.
The coordinate of point $P(v)$ corresponding to some set $\{u,w\}$ with $u,w\in V$ is chosen as follows:
If $v\in\{u,w\}$ and $[u,w]$ is not an edge in $E$, then the coordinate has value~$1$; in all other cases
the coordinate has value $0$.

Since $G$ is $6$-regular, every vertex $v$ has exactly $3q-7$ non-neighbors and hence every point $P(v)$
has exactly $3q-7$ coordinates with value~$1$ (and all other coordinates at~$0$).
Furthermore, if $[u,v]\in E$ then the $L_p$ distance between $P(u)$ and $P(v)$ equals $\ell^*:=\sqrt[p]{6q-14}$,
and if $[u,v]\notin E$ then their $L_p$ distance equals $\sqrt[p]{6q-16}$.
In other words, non-edges correspond to short distances and edges correspond to long distances.
It can be seen that the {\pit} instance has answer YES, if and only if the constructed {\maxtap} instance
has a feasible solution with objective value at least $3q\cdot\ell^*$.
\qed

\begin{lemma}
\label{lp-norm.2}
Problem {\maxtap} with perimeter costs measured according to the Maximum norm $L_{\infty}$ is NP-hard.
\end{lemma}
\proof
The argument is very similar to the argument in Lemma~\ref{lp-norm.1}.
Again we start from an arbitrary instance $G=(V,E)$ of {\pit}, and we create for every vertex $v$ in 
$V_1\cup V_2\cup V_3$ a corresponding point $P(v)$.
We choose the dimension $s=|E|$, and we make every coordinate correspond to one edge in $E$.
For an edge $e=[u,v]\in E$, the coordinates corresponding to $e$ are $0$ for all points with the 
exception of points $P(u)$ and $P(v)$; one of $P(u)$ and $P(v)$ receives coordinate $+1$ and the other 
one receives coordinate $-1$.

Then non-edges correspond to short distances $1$ and edges correspond to long distances~$\ell^*:=2$.
It can be seen that the {\pit} instance has answer YES, if and only if the constructed {\maxtap} 
instance has a feasible solution with objective value at least $6q$.
\qed

\bigskip
\section{A useful lattice}
\label{sec:lattice}
\nopagebreak
In this section, we derive a technical result that will be central in our NP-hardness reduction in
Section~\ref{sec:min}; since this reduction should be implementable on a standard Turing machine,
we want to have all involved numbers to be rational or integer (so that they can be represented by
simple finite strings).
Throughout this section we consider a fixed norm $L_R$ with a fixed unit ball $R$ in the 
Cartesian plane $\rz^2$.

\begin{theorem}
\label{th:lattice}
For any norm $L_R$ with unit ball $R$ in the Cartesian plane $\rz^2$, there exist two integer vectors 
$v_1$ and $v_2$, such that the lattice generated by $v_1$ and $v_2$ has the following properties.
\begin{itemize}
\item[(i)]
The fundamental triangle of the lattice with vertices in $0$, in $v_1$ and in $v_2$ has 
a certain perimeter $\Delta$ (measured in the $L_R$ norm).
\item[(ii)]
Any three (distinct) points $q_1,q_2,q_3$ in the lattice either form a fundamental triangle, or 
otherwise form a triangle with perimeter at least $\Delta+1$ (measured in the $L_R$ norm).
\end{itemize}
\end{theorem}

\begin{figure}[bht]
\bigskip\bigskip
\centerline{\includegraphics[height=3.5cm]{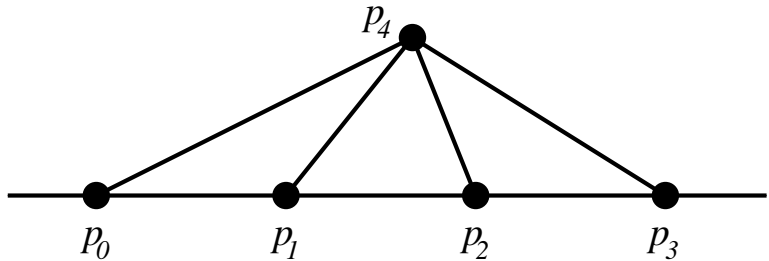}}
\caption{The five points $p_0,p_1,p_2,p_3,p_4$ with the fundamental triangle $p_1p_2p_4$.}
\label{fig:5points}
\end{figure}

The rest of this section is entirely devoted to the proof of Theorem~\ref{th:lattice}.
We start by introducing five points $p_0,p_1,p_2,p_3,p_4$ whose definition is based on a positive 
integer $\aaa$; the value of $\aaa$ will be fixed after the proof of Lemma~\ref{le:lattice.1}.
The Cartesian coordinates of the first four points are given by 
$p_0=(    0,0)$,
$p_1=( \aaa,0)$,
$p_2=(2\aaa,0)$, and
$p_3=(3\aaa,0)$.
These points lie on the $x$-axis, and we assume without loss of generality that
the $L_R$ distance between them is given by 
$d_R(p_i,p_j)=(i-j)\,\aaa$ for all $i,j$ with $0\le i\le j\le 3$.
The final point $p_4$ is chosen in the region above the $x$-axis so that its distances from
$p_0,p_1,p_2,p_3$ satisfy the following inequalities:
\begin{equation}
\label{eq:lattice.0}
\aaa ~<~ d_R(p_1,p_4),~d_R(p_2,p_4) ~\le~ \frac43\aaa ~<~ d_R(p_0,p_4),~ d_R(p_3,p_4)
\end{equation}
See Figure~\ref{fig:5points} for an illustration.

\begin{lemma}
\label{le:lattice.1}
For any integer $\aaa>0$, there exists a point $p_4$ that satisfies the inequalities in 
(\ref{eq:lattice.0}) and that furthermore has rational coordinates.
\end{lemma}
\proof
We let $S$ denote the set of all points $s$ in the upper halfplane that satisfy 
$d_R(p_1,s)=d_R(p_2,s)=4\aaa/3$.
Then set $S$ is the intersection of the boundary of two copies of the unit ball $R$ that are
scaled by the factor $4/3$ and that are centered in points $p_1$ and $p_2$, respectively.
Since $R$ is convex and compact, set $S$ either consists of a single point or otherwise is 
a horizontal line segment.
We claim that $S$ contains some point $s^*$ that simultaneously satisfies 
\begin{equation}
\label{eq:lattice.aux.1}
d_R(p_0,s^*) ~>~ 4\aaa/3 \mbox{\qquad and\qquad}
d_R(p_3,s^*) ~>~ 4\aaa/3.
\end{equation}
First consider the case where $S$ consists of a single point $s=(\bbb,\ccc)$.
Then the horizontal line $\ell$ through this point $s$ contains two points that are at $L_R$ 
distance $4\aaa/3$ from point $p_1$: the point $s$ and some other point that is farther to 
the left of $s$.
(In a degenerate case, the line~$\ell$ contains an entire interval of points whose $L_R$ distance
to $p_1$ equals $4\aaa/3$; in this case point $s$ forms the right endpoint of the interval.) 
If we traverse the points on line $\ell$ from left to right, their distances to point $p_1$ will
follow a convex function; in particular for every point strictly to the right of $s$ the $L_R$ 
distance to $p_1$ will be strictly larger than $d_R(p_1,s)$.
Since the auxiliary point $s'=(\bbb+\aaa,\ccc)$ lies strictly to the right of $s$ on line $\ell$,
we conclude 
\[ d_R(p_1,s') ~>~ d_R(p_1,s) ~=~ 4\aaa/3. \]
Since the line segment $p_0s$ results by shifting line segment $p_1s'$ a distance $\aaa$ to the left,
we derive the desired inequality $d_R(p_0,s)>4\aaa/3$.
A symmetric argument yields $d_R(p_3,s)>4\aaa/3$.
Summarizing, the point $s^*=s$ satisfies the inequalities in (\ref{eq:lattice.aux.1}).

Next consider the case where $S$ is a horizontal line segment between a left endpoint $s_1$ and a
right endpoint $s_2$.
Then the horizontal line $\ell$ through $S$ contains an interval of points whose $L_R$ distance
to $p_1$ equals $4\aaa/3$, and another interval of points whose $L_R$ distance to $p_2$ equals 
$4\aaa/3$; the line segment $S$ is the intersection of these two intervals.
The arguments in the preceding paragraph yield the two inequalities
\begin{equation}
\label{eq:lattice.aux.2}
d_R(p_0,s_2) ~>~ 4\aaa/3 \mbox{\qquad and\qquad} 
d_R(p_3,s_1) ~>~ 4\aaa/3.
\end{equation}
Now let $S_0$ denote the set of all points $s\in S$ with $d_R(p_0,s)\le4\aaa/3$, and let
$S_3$ denote the set of all points $s\in S$ with $d_R(p_3,s)\le4\aaa/3$.
The convexity and the compactness of the unit ball $R$ imply that $S_0$ and $S_3$ are closed intervals.
Furthermore (\ref{eq:lattice.aux.2}) implies $S_0\ne S$ and $S_3\ne S$.
Now suppose for the sake of contradiction that $S=S_0\cup S_3$.
Then the intersection $S_0\cap S_3$ is non-empty and contains a point $t$.
But then the triangle $p_0p_3t$ has one side $p_0p_3$ of length $3\aaa$ and two sides of
length at most $4\aaa/3$.
This is the desired contradiction to the triangle inequality.
We conclude that $S$ contains a point $s^*$ that is neither in $S_0$ nor in $S_3$, and
this point $s^*$ by definition satisfies the desired inequalities (\ref{eq:lattice.aux.1}).

To summarize, we have found a point $s^*\in S$ that satisfies $d_R(p_1,s)=d_R(p_2,s)=4\aaa/3$
and (\ref{eq:lattice.aux.1}).
If $s^*$ has rational coordinates, we are done.
Otherwise, we consider a sufficiently small open neighborhood $N(s^*)$ of $s^*$ whose 
points satisfy (\ref{eq:lattice.aux.1}).
Then the intersection of $N(s^*)$ with the halfplane below $S$ has non-empty interior,
and we can find the desired point with rational coordinates in it.
\qed

\bigskip
Our lattice will have the triangle $p_1p_2p_4$ as fundamental triangle. 
Without loss of generality we assume from now on that the sides of this triangle satisfy
\begin{equation}
\label{eq:lattice.1a}
\aaa ~=~ d_R(p_1,p_2) ~<~ d_R(p_1,p_4) ~\le~ d_R(p_2,p_4).
\end{equation}
Indeed, the first inequality follows from $d_R(p_1,p_2)=\aaa$ and (\ref{eq:lattice.0}),
while the second inequality may be assumed by symmetry.
We now fix the value of $\aaa$ so that $p_4$ has integer coordinates, and so that 
\begin{equation}
\label{eq:lattice.1b}
d_R(p_2,p_4)+1 ~\le~ \min\left\{ d_R(p_0,p_4),~ d_R(p_3,p_4),~ 2\aaa \right\}
\end{equation}
and that
\begin{equation}
\label{eq:lattice.1c}
\aaa+1  ~\le~ d_R(p_2,p_4). 
\end{equation}
The first one of these conditions can be reached by making $\aaa$ a multiple of the denominators of the 
$x$-coordinate and $y$-coordinate of point $p_4$ (which are rational by Lemma~\ref{le:lattice.2}). 
The other conditions can be reached by choosing $\aaa$ sufficiently large so that (\ref{eq:lattice.1b}) 
and (\ref{eq:lattice.1c}) are implied by (\ref{eq:lattice.0}).  
In particular, we will assume from now on that $\aaa\ge3$.

\begin{figure}[bht]
\bigskip\bigskip
\centerline{\includegraphics[height=5.6cm]{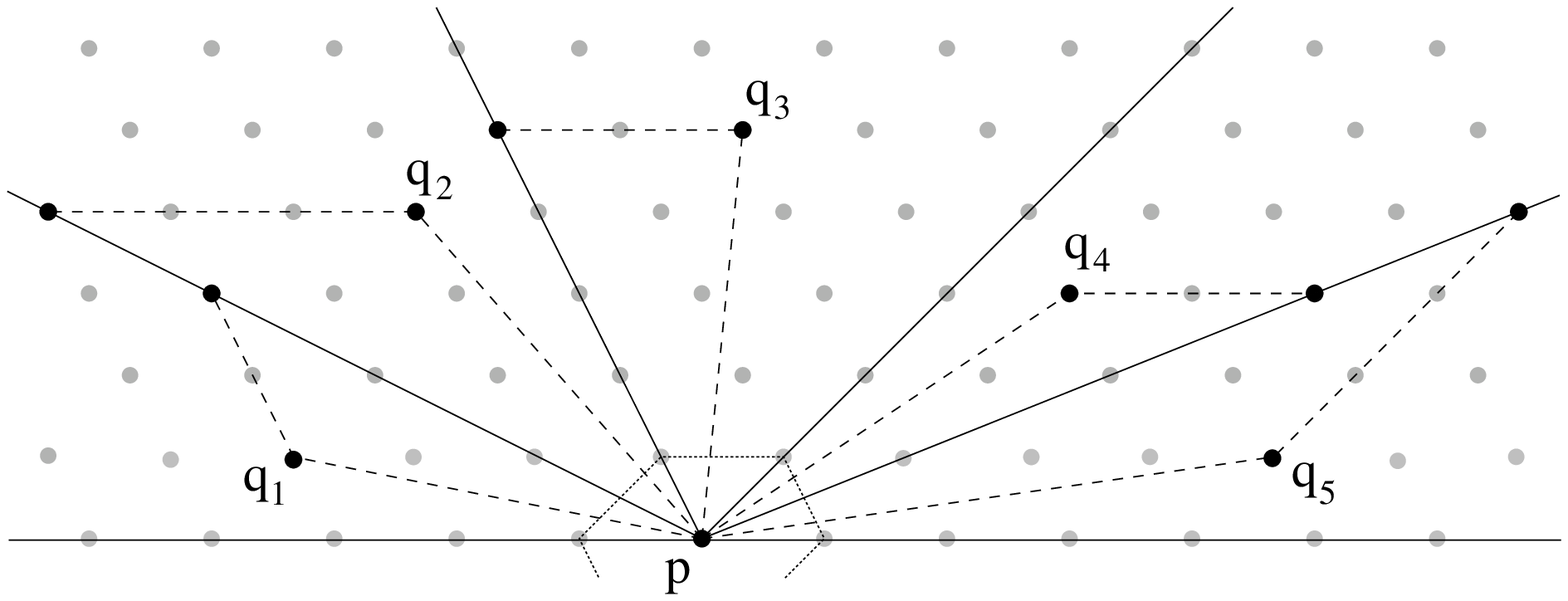}}
\caption{An illustration for the five cases in the proof of Lemma~\protect{\ref{le:lattice.2}}.}
\label{fig:lemma2}
\end{figure}

\begin{lemma}
\label{le:lattice.2}
Let $p$ and $q$ be two points in the lattice with fundamental triangle $p_1p_2p_4$.
If $p$ and $q$ do not both belong to the same fundamental triangle, then $d_R(p,q)\ge d_R(p_2,p_4)+1$.
\end{lemma}
\proof
Without loss of generality we assume that the $y$-coordinate of point $q$ is at least as large as
the $y$-coordinate of point $p$.
We consider the horizontal line through point $p$ together with the four lines through point $p$ 
that are, respectively, parallel to the four line segments 
$p_0p_4$,
$p_1p_4$,
$p_2p_4$, and
$p_3p_4$.
These five lines partition the region above point $p$ into five wedges; 
see Figure~\ref{fig:lemma2} for an illustration.

Let us first deal with the easy cases where point $q$ lies on one of the five lines.
If $q$ lies on the horizontal line then $d_R(p,q)\ge2\aaa$;
if $q$ lies on the line parallel to $p_1p_4$ then $d_R(p,q)\ge2d_R(p_1,p_4)$; and
if $q$ lies on the line parallel to $p_2p_4$ then $d_R(p,q)\ge2d_R(p_2,p_4)$.
In each of these three cases, the desired inequality follows from (\ref{eq:lattice.1a}) and (\ref{eq:lattice.1b}).
Similarly, if $q$ lies on the line parallel to $p_0p_4$ then $d_R(p,q)\ge d_R(p_0,p_4)$, and 
if $q$ lies on the line parallel to $p_3p_4$ then $d_R(p,q)\ge d_R(p_3,p_4)$.
In these two cases the desired inequality follows directly from (\ref{eq:lattice.1b}).

In the main part of the proof we distinguish five cases where point $q$ lies in the interior of one of
the five wedges. 
In the first case, assume that point $q$ lies in the leftmost wedge (like point $q_1$ in 
Figure~\ref{fig:lemma2}).
Draw a line through $q=q_1$ that is parallel to $p_2p_4$, and consider its intersection point $r$
with the upper bounding line of the wedge.
In the triangle $pqr$, the side length $d_R(q,r)$ equals $\lambda$ times $d_R(p_2,p_4)$ and the side length 
$d_R(p,r)$ equals $\mu$ times $d_R(p_3,p_4)$ where $\lambda$ and $\mu$ are positive integers with $\lambda<\mu$.
Now the triangle inequality together with (\ref{eq:lattice.1b}) yields
\begin{eqnarray*}
d_R(p,q) 
&\ge& d_R(p,r) - d_R(q,r) ~~=~~ \mu\cdot d_R(p_3,p_4) - \lambda\cdot d_R(p_2,p_4) \\[0.5ex]
&=&   \lambda\, (d_R(p_3,p_4)-d_R(p_2,p_4)) +(\mu- \lambda)\, d_R(p_3,p_4) \\[0.5ex]
&\ge& d_R(p_3,p_4) ~~\ge~~ d_R(p_2,p_4)+1.
\end{eqnarray*}
This completes the discussion of the first case.
The second, fourth, and fifth case can be handled analogously, and we only list the crucial inequalities for them.
In the second case (where point $q$ lies in the same wedge as point $q_2$ in Figure~\ref{fig:lemma2}), we have
\[ d_R(p,q) ~\ge~ \mu\cdot d_R(p_4,p_3) - \lambda\cdot d_R(p_2,p_3) ~\ge~ d_R(p_4,p_3). \]
In the fourth case (where point $q$ lies in the same wedge as point $q_4$ in Figure~\ref{fig:lemma2}),
\[ d_R(p,q) ~\ge~ \mu\cdot d_R(p_4,p_0) - \lambda\cdot d_R(p_1,p_0) ~\ge~ d_R(p_4,p_0). \]
In the fifth case (where point $q$ lies in the same wedge as point $q_5$ in Figure~\ref{fig:lemma2}),
\[ d_R(p,q) ~\ge~ \mu\cdot d_R(p_0,p_4) - \lambda\cdot d_R(p_1,p_4) ~\ge~ d_R(p_0,p_4). \]
In each of the above three cases, (\ref{eq:lattice.1b}) leads to the desired inequality.
It remains to consider the third case (where point $q$ lies in the same wedge as point $q_3$ in Figure~\ref{fig:lemma2}).
In this case we derive 
\begin{eqnarray*}
d_R(p,q) 
&\ge& \mu\cdot d_R(p_4,p_2) - \lambda\cdot d_R(p_3,p_2) \\[0.5ex]
&=&   \lambda\, (d_R(p_4,p_2)-d_R(p_3,p_2)) +(\mu- \lambda)\, d_R(p_4,p_2).
\end{eqnarray*}
Now $\lambda\ge1$ and $\mu-\lambda\ge1$ together with $d_R(p_4,p_2)-d_R(p_3,p_2)\ge1$ in (\ref{eq:lattice.1c})
yield the desired inequality.
As all five cases have been settled, the proof of the lemma is complete.
\qed

\bigskip
Now let us wrap things up.
Let $\Delta$ denote the $L_R$ perimeter of the fundamental triangle $p_1p_2p_4$.
Note that (\ref{eq:lattice.0}) and (\ref{eq:lattice.1a}) imply the bounds $3\aaa<\Delta\le11\aaa/3$.
By Lemma~\ref{le:lattice.2} and by (\ref{eq:lattice.0}), the three shortest distances between (distinct)
lattice points are the three side lengths $d_R(p_1,p_2)$ and $d_R(p_1,p_4)$ and $d_R(p_2,p_4)$ of the 
fundamental triangle.
All other distances are at least $d_R(p_2,p_4)+1$, that is, the longest side of the fundamental triangle
plus~$1$.

Suppose for the sake of contradiction that for some non-fundamental triangle $r_1r_2r_3$ the $L_R$ perimeter 
would be strictly smaller than $\Delta+1$.
Since this triangle is non-fundamental, by Lemma~\ref{le:lattice.2} one of its side lengths is at least 
$d_R(p_2,p_4)+1$. 
Hence by (\ref{eq:lattice.1a}) and by the above discussion, its two other side lengths must both be equal
to $d_R(p_1,p_2)=\aaa$.
But then the triangle $r_1r_2r_3$ is necessarily degenerate, with all three points on a line and with
$L_R$ perimeter $\aaa+\aaa+2\aaa=4\aaa$.
Now $\aaa\ge3$ implies the desired contradiction $4\aaa\ge 11\aaa/3+1\ge\Delta+1$.
This finally completes the proof of Theorem~\ref{th:lattice}.

\bigskip
\section{The minimization problem}
\label{sec:min}
\nopagebreak
Throughout this section, we investigate versions of {\tap} with perimeter costs where the elements of 
$X\cup Y\cup Z$ are points in the $2$-dimensional Cartesian plane $\rz^2$.
The distances between points are measured according to some fixed norm $L_R$ with unit ball $R$.

\begin{figure}[bht]
\bigskip\bigskip
\centerline{\includegraphics[height=4.5cm]{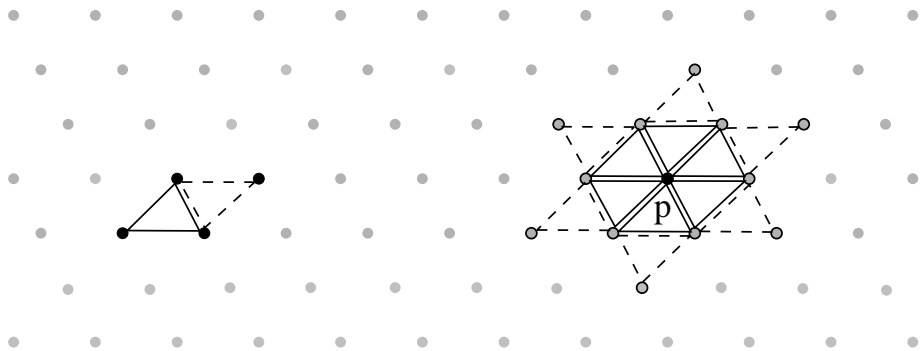}}
\caption{A diamond (to the left) and all possible six directions of a diamond incident to point $p$ (to the right)}
\label{fig:diamond}
\end{figure}

\begin{figure}[bht]
\bigskip\bigskip
\centerline{\includegraphics[height=4.5cm]{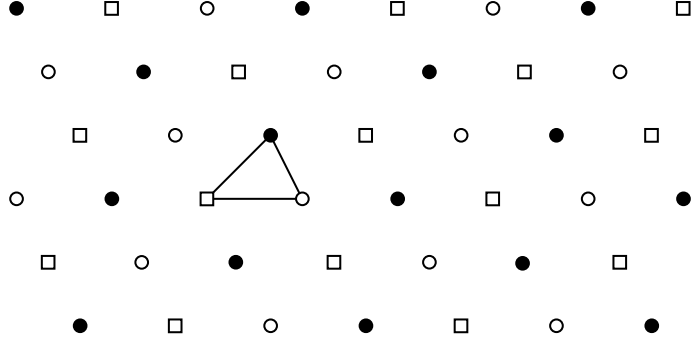}}
\caption{The three-colored lattice}
\label{fig:lattice2}
\end{figure}

We will show that for every (compact, convex, centrally symmetric) unit ball $R$, the resulting version 
of {\tap} with perimeter costs is NP-hard.  
Our reduction is built around the fundamental triangle and the lattice introduced in Theorem~\ref{th:lattice}.
We recall that in this lattice only fundamental triangles have a cheap perimeter of $\Delta$, whereas 
all non-fundamental triangles have an expensive perimeter of at least $\Delta+1$.
A \emph{diamond} is a set of four lattice points obtained by gluing together two fundamental triangles 
along one side; see Figure~\ref{fig:diamond}.
We partition the lattice points into three classes, so that every fundamental triangle contains 
exactly one point from each class.
In the figures the three classes are depicted by circles ($\Circle$), squares ($\Square$) and filled 
circles ($\CIRCLE$); see Figure~\ref{fig:lattice2}.
We refer to this structure as \emph{three-colored lattice}.

Our reduction uses ideas that are similar to those used by
Spieksma \& Woeginger \cite{SpWo1996} and
Pferschy, Rudolf \& Woeginger \cite{PfRuWo1994}.
The reduction is from the following special case of {\tap} whose NP-hardness has been established by
Dyer \& Frieze~\cite{DyFr1986}.
To avoid notational collisions between the variables in {\tap} and the variables in {\ptap}, we 
will consistently denote objects in {\ptap} instances by primed variables.

\begin{quote}
Problem: Planar 3-dimensional assignment problem ({\ptap})
\\[1.0ex]
Instance:
Three pairwise disjoint sets $X'$, $Y'$ and $Z'$ with $|X'|=|Y'|=|Z'|=q'$ and a
set $T'\subseteq X'\times Y'\times Z'$ such that
(i) every element of $X'\cup Y'\cup Z'$ occurs in two or three triples from $T'$, and
(ii) the corresponding graph $G'$ is planar.
(This graph $G'$ contains a vertex for every element of $X'\cup Y'\cup Z'$ and a vertex
for every triple in $T'$. There is an edge connecting a triple vertex to an element vertex
if and only if the corresponding element is a member of the corresponding triple.)
\\[1.0ex]
Question:
Does there exist a subset $T^*$ of $q'$ triples in $T'$ such that each element of
$X'\cup Y'\cup Z'$ is contained in precisely one triple from $T^*$?
\end{quote}

Hence let us consider an arbitrary instance of {\ptap}.
In the first step, we compute a planar layout of the planar graph $G'$ that maps the vertices of $G'$
into integer points in $\zz^2$ and that maps its edges into straight line segments.
This can be done in polynomial time, for instance by using the algorithm of Schnyder \cite{Schnyder1990}.

In the second step, we map the planar layout into the three-colored lattice.
Every point $(\alpha,\beta)$ in the planar layout maps into a point that is in the close neighborhood
of the point $100\alpha\,v_1+100\beta\,v_2$ in the three-colored lattice; here $v_1$ and $v_2$ are the 
integer vectors from Theorem~\ref{th:lattice} that generate the lattice.
Every element of $X'\cup Y'\cup Z'$ is mapped into a corresponding \emph{element point};
every element of $X'$ goes into a circle ($\Circle$), 
every element of $Y'$ goes into a square ($\Square$), and 
every element of $Z'$ goes into a filled circle ($\CIRCLE$).
Every triple in $T'$ is mapped into a fundamental triangle called \emph{triple triangle}.
These element points and triple triangles roughly imitate the planar layout constructed above;
there is plenty of leeway for doing this, since the main restriction is that the various objects 
should be embedded far away from each other.

\begin{figure}[bht]
\bigskip\bigskip
\centerline{\includegraphics[height=4.2cm]{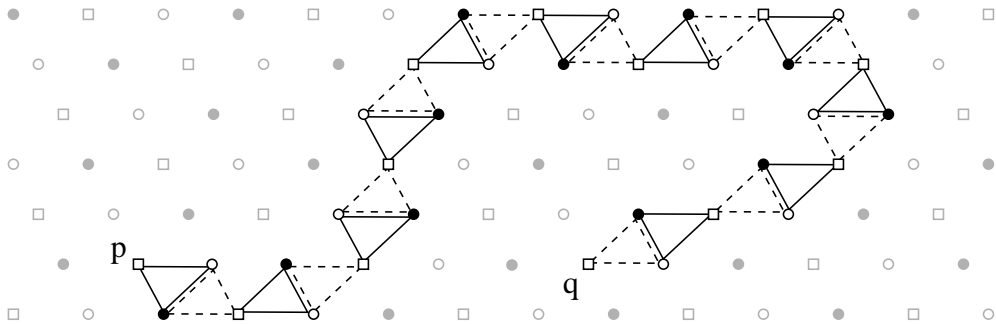}}
\caption{A chain of diamonds between two points $p$ and $q$}
\label{fig:chain}
\end{figure}

\begin{figure}[bht]
\bigskip\bigskip
\centerline{\includegraphics[height=5.2cm]{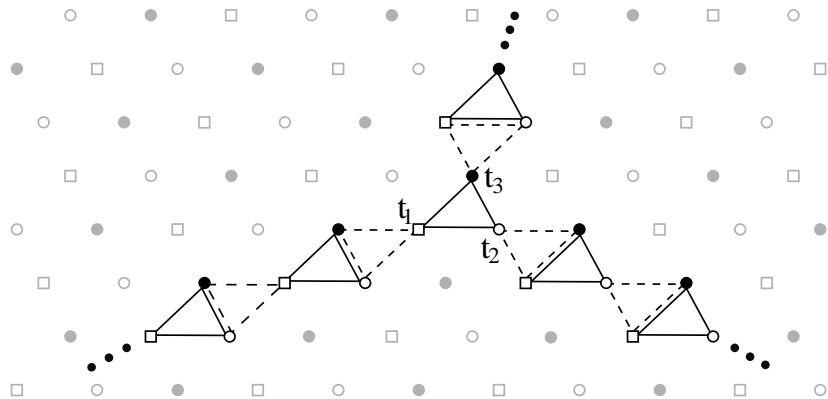}}
\caption{How chains of diamonds attach to a triple triangle}
\label{fig:trGadget}
\end{figure}

\begin{figure}[bht]
\bigskip\bigskip
\centerline{\includegraphics[height=5.2cm]{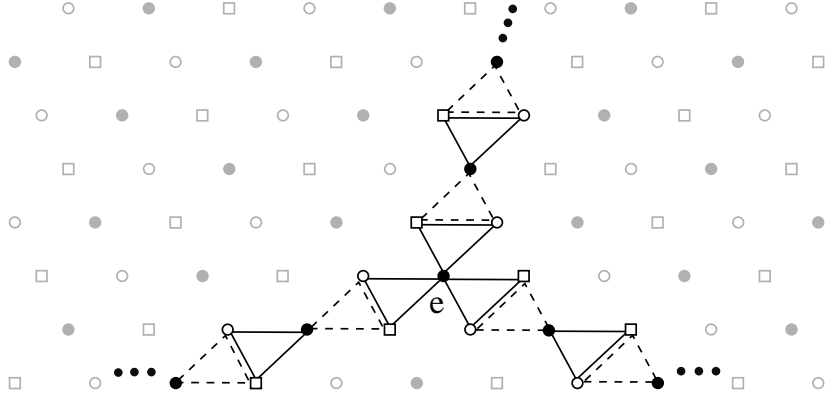}}
\caption{How chains of diamonds attach to an element point}
\label{fig:elGadget}
\end{figure}

In the third step, we introduce several chains of diamonds that connect certain element points
to certain triple triangles; see Figure~\ref{fig:chain} for an illustration.
Every such chain connects an element point (for some element $x'$ of $X'\cup Y'\cup Z'$) to a triple
triangle (whose corresponding triple $t'$ in $T'$ contains that element $x'$).
These chains roughly follow the straight line segment that corresponds to the edge between $x'$ and
$t'$ in the planar layout in the first step.
Figure~\ref{fig:trGadget} shows how such a chain is attached to a triple triangle, and
Figure~\ref{fig:elGadget} shows how such a chain is attached to an element point.

\begin{figure}[bht]
\bigskip\bigskip
\centerline{\includegraphics[height=3.8cm]{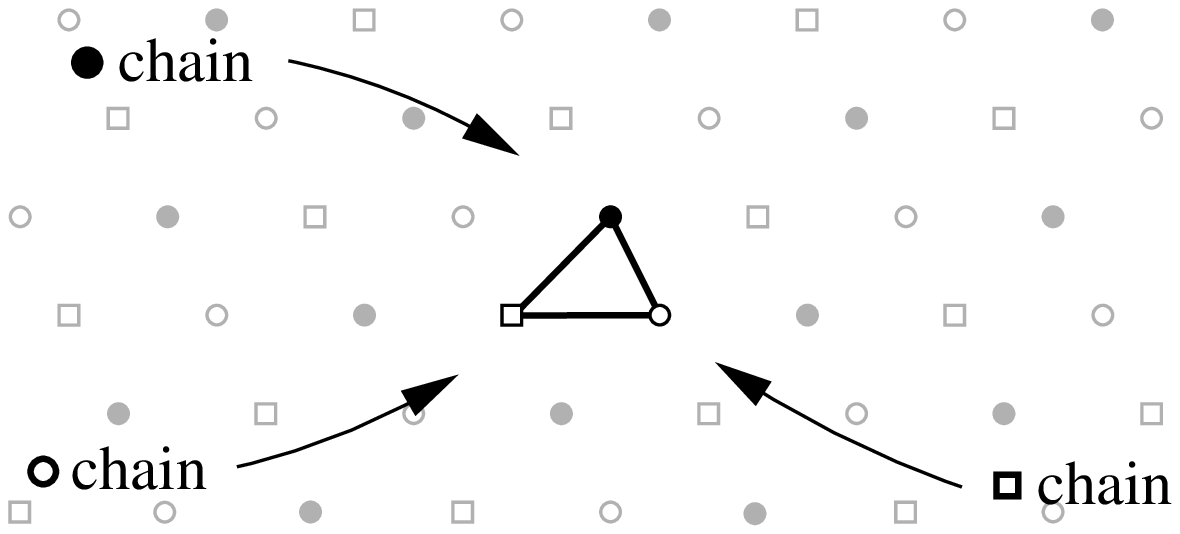}}
\bigskip\bigskip\bigskip
\centerline{\includegraphics[height=3.8cm]{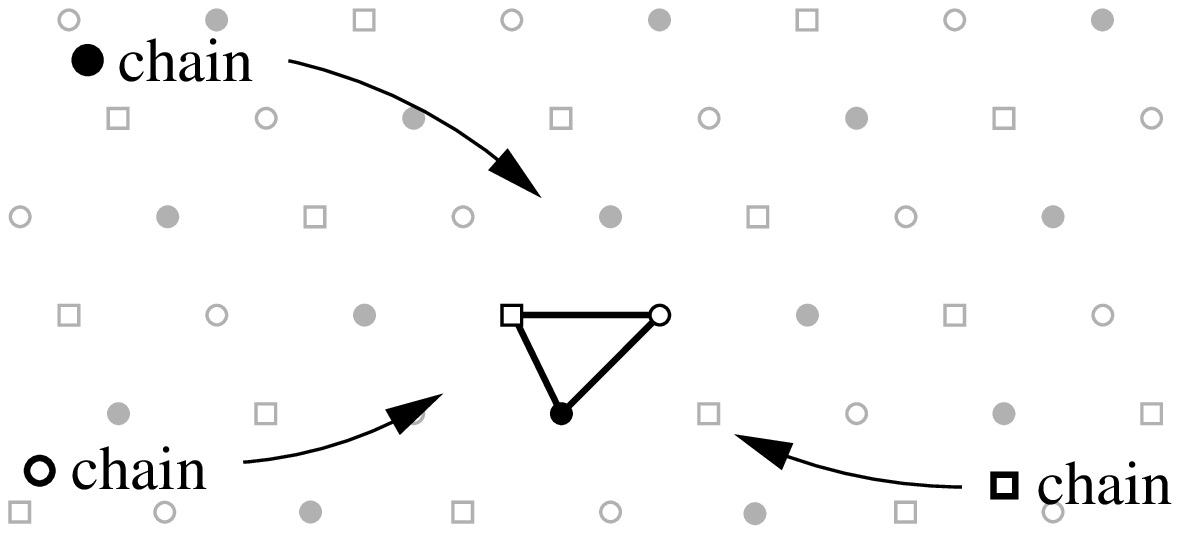}}
\caption{Connecting chains to triple triangles: the upper picture shows an infeasible clockwise choice,
the picture at the bottom shows the feasible counter-clockwise choice}
\label{fig:orientations}
\end{figure}

Three comments are in place.
First, if an element $x'$ of $X'\cup Y'\cup Z'$ occurs in only two triples in $T'$, then the
corresponding element point is attached to only two chains of diamonds.
Secondly, for every chain of diamonds the attachment point in the triple triangle belongs to the same
class ($\Circle$, $\Square$, $\CIRCLE$) as the element point at the other end of the chain.
Thirdly, we note that there are two combinatorially different ways of choosing a triple triangle in the 
lattice; one way has the vertices in the classes $\Circle$, $\Square$, $\CIRCLE$ clockwise, and the other 
way has the vertices in the classes $\Circle$, $\Square$, $\CIRCLE$ counter-clockwise.
We always pick the way that allows a crossing-free attachment of the three
chains of diamonds 
to the triple triangle; see Figure~\ref{fig:orientations} for an illustration.

The element points, triple triangles and chains of diamonds altogether contain $3n$ points from 
the three-colored lattice, and each of the three classes contains exactly $n$ points.
These three sets with $n$ points form the three sets $X$, $Y$, $Z$ in a {\tap} instance with
perimeter costs.
We complete the reduction by defining the integer bound $B=\lceil n\,\Delta\rceil$.
The following two lemmas establish the connections between the considered instance of {\ptap}
and the newly constructed instance of {\tap}.

\begin{lemma}
\label{le:min.1}
If the constructed instance of {\tap} has a solution with objective value at most $B$,
then the considered instance of {\ptap} has answer YES.
\end{lemma}
\proof
Assume that the {\tap} instance has a solution with objective value at most $B$.
Then by Theorem~\ref{th:lattice} all $n$ triples in this solution have perimeter cost $\Delta$
and induce fundamental triangles in the lattice.
Moreover, it is straightforward to verify that from any chain of diamonds the solution does 
either pick all the dashed triangles or does pick all the solid triangles; 
see Figures~\ref{fig:trGadget} and~\ref{fig:elGadget}.

We define a subset $T^*$ of the triples in $T'$ by picking all the triples for which the
corresponding triple triangle occurs in the solution for the {\tap} instance.
Consider some element $x'\in X'\cup Y'\cup Z'$.
The corresponding element point is contained in exactly one solid triangle in the {\tap} solution,
and this triangle must belong to some chain; see Figure~\ref{fig:elGadget}. 
Consequently, this is a chain of solid triangles which propagates to some triple triangle. 
Figure~\ref{fig:trGadget} shows that the corresponding triple triangle is in $T^*$.
To summarize, every element $x'\in X'\cup Y'\cup Z'$ is contained in exactly one triple in $T^*$.
Hence $T^*$ yields the desired certificate that the {\ptap} instance has answer YES.
\qed

\begin{lemma}
\label{le:min.2}
If the considered instance of {\ptap} has answer YES, then the constructed instance of {\tap} has
a solution with objective value at most $B$.
\end{lemma}
\proof
Assume that the {\ptap} instance has answer YES, so that there is a set $T^*$ of $q'$ triples in $T'$
that covers every element of $X'\cup Y'\cup Z'$ exactly once.
Then we construct the following solution for the {\tap} instance.
For every triple in $T^*$, we pick the corresponding triple triangle for the {\tap} solution.
For every element $x'\in X'\cup Y'\cup Z'$, we pick the solid triangles in the chain of diamonds
that connects the element point for $x'$ to the triple triangle whose triple covers $x'$ in $T^*$;
in the other chains incident to this element point, we pick the dashed triangles.
As all points in $X\cup Y\cup Z$ are covered by the $n$ picked (fundamental!) triangles,
their overall length equals $n\,\Delta$.
\qed

\bigskip
Note that the bound $B$ in our construction is integer, and note that all the points in $X\cup Y\cup Z$ 
have integer coordinates; hence the reduction can easily be implemented in polynomial time (and without
worrying about computations with irrational numbers).
Together with Lemmas~\ref{le:min.1} and \ref{le:min.2} this yields the following theorem.

\begin{theorem}
\label{th:min.main}
For any fixed norm $L_R$ with unit ball $R$ in two-dimensional space $\rz^2$, problem {\tap} with perimeter
costs measured according to $L_R$ is NP-hard.
\qed
\end{theorem}

\bigskip
\section{Implications for the 3-dimensional matching problem}
\label{sec:3dm}
\nopagebreak
Up to this point we have been solely concerned with the 3-dimensional \emph{assignment} problem, where
the underlying elements belonged to three classes $X$, $Y$ and $Z$, and where every triple contained
exactly one element from every class.
In the closely related 3-dimensional \emph{matching} problem ({\tdm}) all the elements belong to the same class:
An instance of {\tdm} consists of a ground set $U$ with $|U|=3n$ and a cost function $c:U\times U\times U\to\rz$.
The goal is to find a set of $n$ triples in $U\times U\times U$ that cover every element in $U$ exactly once, 
such that the sum of the costs of these triples is minimized.
In the maximization version {\maxtdm} of {\tdm}, this sum is to be maximized.

The algorithmic behavior of {\tdm} is very similar to that of {\tap}.
Both problems are NP-hard in general, and (as a rule of thumb) algorithms for one problem usually translate into
similar algorithms for the other problem.
Pferschy, Rudolf \& Woeginger \cite{PfRuWo1994} proved that {\tdm} with perimeter costs under Euclidean distances
in $\rz^2$ is NP-hard.
Our hardness arguments in Section~\ref{sec:min} can easily be adapted to {\tdm} by setting $U:=X\cup Y\cup Z$,
thus extending and generalizing the result of \cite{PfRuWo1994} to arbitrary norms.

\begin{corollary}
\label{co:3dm.1}
For any fixed norm $L_R$ with unit ball $R$ in two-dimensional space $\rz^2$, problem {\tdm} with perimeter
costs measured according to $L_R$ is NP-hard.
\qed
\end{corollary}

Also the NP-hardness proofs in Section~\ref{sec:lp-norm} for {\maxtap} (when the dimension is part of the
input) can easily be carried over to the matching problem.
\begin{corollary}
\label{co:3dm.2}
For any fixed $p$ with $1\le p\le\infty$, problem {\maxtdm} with perimeter costs measured according 
to the $L_p$ norm is NP-hard.
\qed
\end{corollary}

In a similar fashion, the positive results in Sections~\ref{sec:tunnel} and \ref{sec:polynorm} for the 
maximization version carry over to the 3-dimensional matching problem. 
We leave the (fairly easy) technical details to the reader.
\begin{corollary}
\label{co:3dm.3}
Problem {\maxtdm} with perimeter costs under tunneling distances can be solved within a time complexity that
depends polynomially on the instance size $n$ (and exponentially on the number of tunnels).
\qed
\end{corollary}

\begin{corollary}
\label{co:3dm.4}
For any polyhedral norm $L_R$ with unit ball $R$ in $s$-dimensional space $\rz^s$, problem {\maxtdm} with
perimeter costs measured according to $L_R$ can be solved within a time complexity that depends polynomially
on the instance size $n$ (and exponentially on the number of facets of the polyhedral unit ball).
\qed
\end{corollary}

\begin{corollary}
\label{co:3dm.5}
For any fixed (not necessarily polyhedral) norm $L_R$ with unit ball $R$ in $s$-dimensional space $\rz^s$,
problem {\maxtdm} with perimeter costs measured according to $L_R$ possesses a PTAS.
\qed
\end{corollary}

\bigskip
\section{Conclusions}
\label{sec:conclusions}
\nopagebreak
We have derived a variety of results on the complexity of {\tap} and {\maxtap} with perimeter costs,
when distances are measured according to certain norms.

Problem {\tap} turned out to be hard for all norms, even if the dimension of the underlying 
Cartesian space $\rz^s$ equals $s=2$.
This of course (trivially) implies NP-hardness also for all dimensions $s\ge3$.
Problem {\maxtap} with perimeter costs shows a more versatile behavior.
If the dimension $s$ is fixed then {\maxtap} is easy for polyhedral norms. 
If the dimension $s$ is part of the input then {\maxtap} is NP-hard for any $L_p$ norm.  
The following question does not seem to be within the reach of our methods.
\begin{open}
Decide whether {\maxtap} with perimeter costs is NP-hard, if the elements are points in $2$-dimensional 
space $\rz^2$ and if the distances are measured according to the Euclidean norm $L_2$.
\end{open}

The literature contains only a handful of results on the approximability of {\tap} and {\maxtap} with 
perimeter costs.
Our Theorem~\ref{th:polynorm.ptas} yields the existence of a PTAS for {\maxtap} if the dimension $s$ is fixed.
Furthermore, there is a polynomial time approximation algorithm with worst case guarantee $4/3$ for the 
{\tap} with perimeter costs by Crama \& Spieksma \cite{CrSp1992}, which works for arbitrary metric spaces.
The following open problem seems to be very challenging.

\begin{open}
Establish APX-hardness of the minimization problem {\tap} with perimeter costs, if the elements are points 
in $2$-dimensional space $\rz^2$ and if the distances are measured according to the Euclidean norm $L_2$.
\end{open}

\bigskip
\paragraph{Acknowledgements.}
This research was supported by the Austrian Science Fund: W1230, Doctoral Program ``Discrete Mathematics''.
Gerhard Woeginger acknowledges support
by DIAMANT (a mathematics cluster of the Netherlands Organization for Scientific Research NWO),
and by the Alexander von Humboldt Foundation, Bonn, Germany.



\bigskip
\begin{thebibliography}{18}

\bibitem{Barvinok1996}
{\sc A.I. Barvinok} (1996).
Two algorithmic results for the traveling salesman problem.
\emph{Mathematics of Operations Research 21}, 65--84.

\bibitem{Barvinok2003}
{\sc A.I. Barvinok, S.P. Fekete, D.S. Johnson, A. Tamir, G.J. Woeginger and R. Woodroofe} (2003).
The geometric maximum travelling salesman problem.
\emph{Journal of the ACM 50}, 641--664.

\bibitem{Burkard-book}
{\sc R.E. Burkard, M. Dell'Amico, and S. Martello} (2009).
\emph{Assignment Problems}.
SIAM, Philadelphia.

\bibitem{BuRuWo1996}
{\sc R.E. Burkard, R. Rudolf and G.J. Woeginger} (1996).
Three-dimensional axial assignment problems with decomposable cost coefficients.
\emph{Discrete Applied Mathematics 65}, 123--140.

\bibitem{CrSp1992}
{\sc Y. Crama and F.C.R. Spieksma} (1992).
Approximation algorithms for three-dimensional assignment problems with triangle inequalities.
\emph{European Journal of Operational Research 60}, 273--379.

\bibitem{DyFr1986}
{\sc M.E. Dyer and A.M. Frieze} (1986).
Planar 3DM is NP-complete.
\emph{Journal of Algorithms 7}, 174--184.

\bibitem{GaJo1979}
{\sc M.R. Garey and D.S. Johnson} (1979).
\emph{Computers and Intractability: A Guide to the Theory of NP-Completeness}.
Freeman, San Francisco.

\bibitem{Karp1972}
{\sc R.M. Karp} (1972).
Reducibility among combinatorial problems.
In R.E. Miller and J.W. Thatcher (eds), \emph{``Complexity of Computer Computations''}, 
Plenum Press, New York, 85--103.

\bibitem{Lenstra1983}
{\sc H.W. Lenstra} (1983).
Integer programming with a fixed number of variables.
\emph{Mathematics of Operations Research 8}, 538--548.

\bibitem{PfRuWo1994}
{\sc U. Pferschy, R. Rudolf and G.J. Woeginger} (1994).
Some geometric clustering problems.
\emph{Nordic Journal of Computing 1}, 246--263.

\bibitem{PoSpWo2013}
{\sc S. Polyakovskiy, F.C.R. Spieksma and G.J. Woeginger} (2013).
The three-dimensional matching problem in Kalmanson matrices. 
\emph{Journal of Combinatorial Optimization 26}, 1--9.

\bibitem{Schnyder1990}
{\sc W. Schnyder} (1990).
Embedding planar graphs on the grid.
\emph{Proceedings of the 1st Annual ACM-SIAM Symposium on Discrete Algorithms (SODA'1990)}, 138--147.

\bibitem{SpWo1996}
{\sc F.C.R. Spieksma and G.J. Woeginger} (1996).
Geometric three-dimensional assignment problems.
\emph{European Journal of Operational Research 91}, 611--618.

\bibitem{VanRooij2013}
{\sc J.M.M. van Rooij, M.E. van Kooten Niekerk and H.L. Bodlaender} (2013).
Partition into triangles on bounded degree graphs. 
\emph{Theory of Computing Systems 52}, 687--718.

\end{thebibliography}
\end{document}